\documentclass[runningheads]{llncs}
\usepackage{graphicx}
\usepackage{amsmath}
\usepackage{xcolor}
\begin{document}
\title{MeshingNet: A New Mesh Generation Method based on Deep Learning}
\author{Zheyan Zhang\orcidID{0000-0002-8779-6659} \and
Yongxing Wang\orcidID{0000-0002-5673-042X}
\and\\
Peter K. Jimack\orcidID{0000-0001-9463-7595} \and 
He Wang\orcidID{0000-0002-2281-5679} }

\authorrunning{Z. Zhang et al.}
\institute{School of Computing, University of Leeds, Leeds, UK, LS2 9JT}
\maketitle              % typeset the header of the contribution
\begin{abstract}
We introduce a novel approach to automatic unstructured mesh generation using machine learning to predict an optimal finite element mesh for a previously unseen problem. The framework that we have developed is based around training an artificial neural network (ANN) to guide standard mesh generation software, based upon a prediction of the required local mesh density throughout the domain. We describe the training regime that is proposed, based upon the use of \emph{a posteriori} error estimation, and discuss the topologies of the ANNs that we have considered. We then illustrate performance using two standard test problems, a single elliptic partial differential equation (PDE) and a system of PDEs associated with linear elasticity. We demonstrate the effective generation of high quality meshes for arbitrary polygonal geometries and a range of material parameters, using a variety of user-selected error norms.

\end{abstract}

\section{Introduction}
\label{Introduction}
Mesh generation is a critical step in the numerical solution of a wide range of problems arising in computational science. The use of unstructured meshes is especially common in domains such as computational fluid dynamics (CFD) and computational mechanics, but also arises in the application of finite element (FE) and finite volume (FV) methods for estimating the solutions of general partial differential equations (PDEs): especially on domains with complex geometries~\cite{shewchuk2002delaunay,si2015tetgen}. The quality of the FE/FV solution depends critically on the nature of the underlying mesh. For an ideal mesh the error in the FE solution (we focus on the FE method (FEM) in this paper) will be distributed equally across the elements of the mesh, implying the need for larger elements where the local error density is small and smaller elements where the local error density is large. This “equidistribution” property tends to ensure that a prescribed global error tolerance (i.e.\ an acceptable error between the (unknown) true solution and the computed FE solution) can be obtained with the fewest number of elements in the mesh~\cite{dorfler1996convergent,stevenson2007optimality}. This is a desirable feature since the computational work generally grows superlinearly with the number of elements (though, in some special cases, this can be linear~\cite{notay2010aggregation}).

The conventional approach to obtain high quality meshes involves multiple passes, where a solution is initially computed on a relatively coarse uniform mesh and then a post-processing step, known as \emph{a posteriori} error estimation, is undertaken~\cite{ainsworth1997posteriori,bank1985some,zienkiewicz1991adaptivity}. This typically involves solving many auxiliary problems (e.g.\ one per element or per patch of elements) in order to estimate the local error in the initial solution \cite{gratsch2005posteriori}. These local errors can be combined to form an overall (global) error estimate but they can also be used to determine where the local mesh density most needs to be increased (mesh refinement), and by how much, and where the local mesh density may be safely decreased (mesh coarsening), and by how much. A new mesh is then generated based upon this \emph{a posteriori} error estimate and a new FE solution is computed on this mesh. A further \emph{a posteriori} error estimate may be computed on this new mesh to determine whether the solution is satisfactory (i.e. has an error less than the prescribed tolerance) or if further mesh refinement is required.

A necessary requirement for efficient \emph{a posteriori} error estimators is that they should be relatively cheap to compute (whilst still, of course, providing reliable information about the error in a computed solution). For example, the approaches of \cite{ainsworth1997posteriori,apel2004,bank1985some} each solve a supplementary local problem on each finite element in order to estimate the 2-norm of the local error from the local residual. Alternatively, recovery-based error estimators use local ``superconvergence'' properties of finite elements to estimate the energy norm of the local error based purely on a locally recovered gradient: for example the so-called ZZ estimator of~\cite{zienkiewicz1991adaptivity}. Here, the difference between the original gradient and a patch-wise recovered gradient indicates the local error. In the context of linear elasticity problems, the elasticity energy density of a computed solution is evaluated at each element and the recovered energy density value at each vertex is defined to be the average of its adjacent elements. The local error is then proportional to the difference between the recovered piece-wise linear energy density and the original piece-wise constant values. 

In this paper we exploit a data-driven method to improve the efficiency of non-uniform mesh generation compared with existing approaches. The core of non-uniform mesh generation is to find an appropriate mesh density distribution in space. Rather than utilizing expensive error estimators at each step, we compute and save high quality mesh density distributions obtained by FEM, followed by accurate error estimation, as a pre-processing step. If a model can successfully learn from the data, it no longer needs an FE solution and error estimator to predict a good mesh density distribution, but instead can reply on learning from a set of similar problems for prediction. Artificial Neural Networks (ANNs) are mathematical models that use a network of ``neurons'' with activation functions to mimic biological neural networks. Even a simple ANN can approximate continuous functions on compact subsets of $R^n$~\cite{csaji2001approximation}.
 An ANN is composed of a large number of free parameters which define the network that connects the ``neurons''. These trainable parameters are generally not explainable. However, with them ANNs can approximate the mapping between inputs and outputs. A training loss function reflects how well the predicted output of an ANN performs for a given input (i.e.\ measured by the difference between the ANN's prediction and the ground truth). Furthermore, an ANN can be trained by gradient decent methods because this loss function is generally differentiable with respect to the network parameters.
In recent years, With the developments of parallel hardware, larger/deeper neural networks (DNNs) have been proven to supersede existing methods on various high-level tasks such as object recognition~\cite{krizhevsky2012imagenet}. Within computational science, DNNs have also been explored to solve ordinary differential equations (ODEs) and PDEs under both supervised~\cite{chen2018neural,long2017pde} and unsupervised~\cite{han2018solving} settings.

In the work reported here we propose a DNN model, MeshingNet, to learn from the \emph{a posteriori} error estimate on an initial (coarse) uniform mesh and predict non-uniform mesh density for refinement, without the need for (or computational expense of) solving an FE system or computing an \emph{a posteriori} error estimate. MeshingNet is trained using an accurate error estimation strategy which can be expensive but is only computed offline. Hence, the mesh generation process itself is extremely fast since it is able to make immediate use of standard, highly-tuned, software (in our case~\cite{shewchuk2002delaunay}) to produce a high quality mesh at the first attempt (at similar cost to generating a uniform mesh with that number of elements). Note that it is not our intention in this work to use deep learning to solve the PDEs directly (as in~\cite{long2017pde,sirignano2018dgm} for example): instead we simply aim to provide a standard FE solver with a high quality mesh, because we can provide more reliable  predictions in this way, based upon the observation that a greater variation in the quality of predictions can be tolerated for the mesh than for the solution itself. For example, in an extreme worst case where the DNN predicts a constant output for all inputs, the result would be a uniform mesh (which is tolerable) however such a poor output would be completely unacceptable in the case where the network is used to predict the solution itself.

Formally, we propose, what is to the best of our knowledge, the first DNN-based predictor of \emph{a posteriori} error, that can: (i) efficiently generate non-uniform meshes with a desired speed; (ii) seamlessly work with existing mesh generators, and; (iii) generalize to different geometric domains with various governing PDEs, boundary conditions (BCs) and parameters.

The remainder of this paper is structured as follows. In the next section we describe our proposed deep learning algorithm. This is not the first time that researchers have attempted to apply ANNs to mesh generation. However,  previous attempts have been for quite specific problems~\cite{chedid1996automatic,lowther1993density} and have therefore been able to assume substantially more \emph{a priori }knowledge than our approach. Consequently, the novelty in section 2 comes through both the generality of the approach used in formulating the problem (i.e.\ generality of the inputs and outputs of the DNN) as well as the network itself. In Section 3, we demonstrate and assess the performance of our approach on two standard elliptic PDE problems. These tests allow us to account for variations in the PDE system, the domain geometry, the BCs, the physical problem parameters and the desired error norm when considering the efficacy of our approach. Finally, in Section 4 we discuss our plans to further develop and apply this deep learning approach.

\section{Proposed method}
\label{Method}
\subsection{Overview}
\label{Overview}
We consider a standard setting where the FEM is employed. Given a geometry and a mesh generator, a low density uniform mesh (LDUM) can be easily computed, then refined non-uniformly based on the \emph{a posteriori} error distribution, for better accuracy. Since this iterative meshing process is very time-consuming, we propose a DNN-based supervised learning method to accelerate it.

\begin{figure}[]
\includegraphics[scale =0.265]{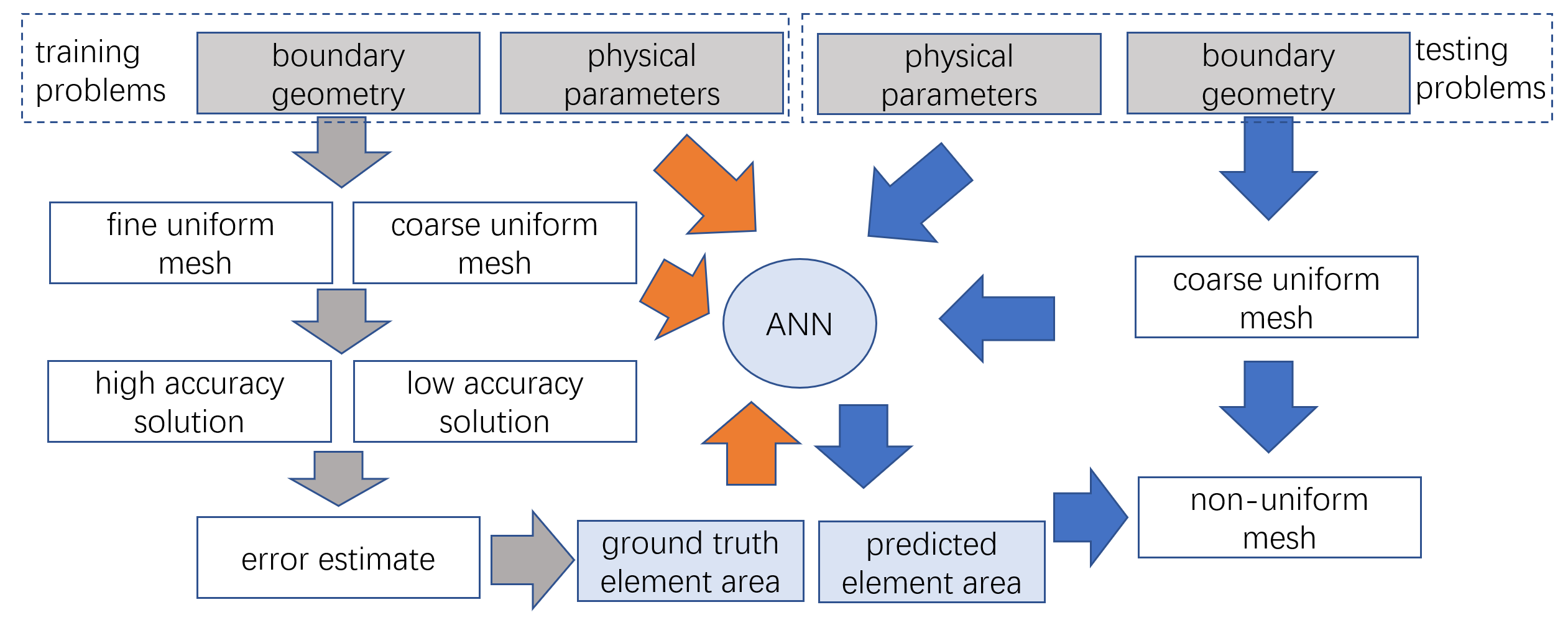}
\caption{A high-level diagram of MeshingNet workflow showing the ANN alongside the training regime (left) and testing regime (right). Grey, orange and blue arrows represent data generation, training and testing processes.}
\label{fig 7}
\end{figure}

Given a family of governing PDEs and material parameters, we assume that there is a mapping
\begin{equation}
F: \Gamma,B,M,X\rightarrow A(X)
\label{mapping}
\end{equation}
that can be learned, where $\Gamma$ is a collection of domain geometries, $B$ is a set of BCs, $M$ is a set of PDE parameters (e.g.\ material properties), $x \in X$ is a location in the domain and $A(X)$ is the target area upper bound distribution over the whole domain. To represent an interior location, we use Mean Value Coordinates \cite{floater2003mean} because they provide translational and rotational invariance with respect to the boundary. Given $\Gamma$, $B$, $M$ and $X$, we aim to predict $A(X)$ quickly. The mapping $F$ is highly non-linear and is therefore learned by our model, MeshingNet. Under the supervised learning scheme, we first build up our training data set by computing high-accuracy solutions (HASs) on high-density uniform meshes (HDUMs) using a standard FE solver. The same computation is also done on LDUMs to obtain lower accuracy solutions (LAS). Then an \emph{a posteriori} error distribution $E(X)$ is computed based upon interpolation between these solutions. According to $E(X)$, we compute $A(X)$ for refinement. The training data is enriched by combining different geometries with different parameters and BCs. Next, MeshingNet is trained, taking as input the geometry, BCs and material properties, with the predicted local area upper bound $A(X)$ as output. After training, MeshingNet is used to predict $A(X)$ over a new geometry with a LDUM. The final mesh is generated either by refining the LDUM non-uniformly or using it to generate a completely new mesh (e.g. using the method in \cite{shewchuk2002delaunay}), guided by the predicted target local area upper bound.
Fig.~\ref{fig 7} illustrates the whole workflow of our mesh generation system.

The approach that we propose has a number of components that are not fixed and may be selected by the user. MeshingNet is agnostic about both the mesh generator and the particular FE/FV solver that are used. It is designed to work with existing methods. Furthermore, the \textit{a posteriori} error can be computed using any user-defined norm (in this paper we consider L1 and energy norms respectively in our two validation tests). Some specific examples of governing equations, geometries, boundary conditions and material parameters are illustrated in the evaluation section. Prior to this however we provide some additional details of the components used in our paper.

\subsection{Component Details}
\label{Components}
\textbf{Mesh generation.}
All meshes (LDUM, HDUM and the refined mesh) are created using the software \textit{Triangle}~\cite{shewchuk2002delaunay} which conducts Delaunay triangulations.  \textit{Triangle} reads a planar graph, representing the boundary of a two-dimensional polygonal domain as its input and the user may define a single maximum element area across the entire polygon when a uniform mesh is desired. To refine a mesh non-uniformly the user specifies, within each coarse element, what the element area upper bound should be after refinement. \textit{Triangle} ensures a smooth spatial variation of the element size. The refinement is not nested since, although it does not eliminate the pre-existing vertices of the original mesh, it does break the original edges.

\textbf{Mesh refinement via error estimation.}
Broadly speaking, the finer the mesh is, the closer the FE result is expected to be to the ground truth. We regard the HAS as the ground truth by ensuring that HDUMs are always significantly finer than LDUMs. Linear interpolation is used to project the LAS to the fine mesh to obtain LAS*, where LAS* has the same dimension as HAS. An error estimate approximates the error $E$ by comparing LAS* and HAS in the selected norm, on the HDUM, and this is then projected back to the original LDUM. The target area upper bound for the refined mesh within each LDUM element is then defined to be inversely correlated with $E$.
In this paper, we select $K/E(x_i)^\alpha$ as the area upper bound for element number $i$ of the LDUM (to be refined), where $x_i$ is the center of the $i$ th element, $K$ and $\alpha$ determines the overall target element number and, in the examples given here, we always choose $\alpha = 1$.
By varying $K$ appropriately it is possible to adjust the refined mesh to reach a target total number of elements.
%Care should be taken however not to choose $K$ to be so small that the area of any elements within the refined mesh is finer than (or even as fine as) the ground truth mesh.

\begin{figure}[]
\includegraphics[scale =0.55]{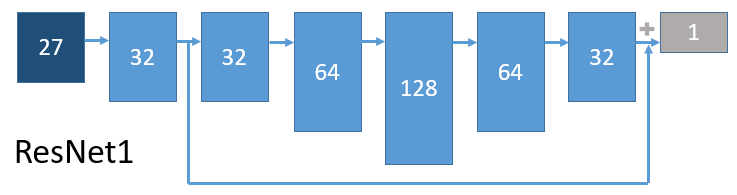}
\includegraphics[scale =0.55]{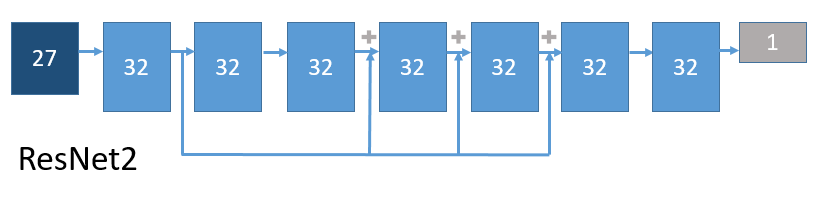}
\caption{Our two residual networks, illustrated with 27 input parameters. ResNet1 is a modification of the FCN with the output of the first hidden layer added to the output of the last hidden layer. ResNet2 has all hidden layers of the same dimension and the output of the first hidden layer is added to outputs of the three front hidden layers.}
\label{fig4}
\end{figure}

\textbf{MeshingNet model and training.}
As outlined in subsection~\ref{Overview}, for a given PDE system, MeshingNet approximates the target local element area upper bound $A(x)$ at a given point within a polygonal domain based upon inputs which include: the coordinates of the polygon's vertices, key parameters of the PDE and the mean value coordinates of the specified point (mean value coordinates, \cite{floater2003mean}, parameterizes a location within a 2D polygon as a convex combination of polygon vertices). Two types of DNNs are considered: a fully connected network (FCN) and two residual networks (ResNets). The dimensions of our FCN layers are X-32-64-128-128-64-32-8-1 (where X represents the dimension of the input and is problem-specific) and each hidden layer uses rectified linear units (ReLU) as the activation function. To further improve and accelerate training, two ResNets are also experimented with to enhance FCN (Fig.~\ref{fig4}). ResNet1 enhances FCN by adding a connection from the first hidden layer to the output of the last one. Note that residual connections can help to resolve the vanishing gradient problem in deep networks and improve training \cite{veit2016residual}. ResNet2 enhances FCN by adding multiple residual connections. This is inspired by recent densely connected convolutional networks \cite{Huang2017dense} which shows superior data-fitting capacity with a relatively small number of parameters. The training data set samples over geometries, BCs and parameter values. Each geometry with fixed BCs and parameters uniquely defines a problem. In a problem, each LDUM element centroid (represented by its mean value coordinates) and its target $A(x)$ forms an input-output training pair. We randomly generated 3800 problems of which 3000 are used for training and 800 for testing (because each LDUM contains approximately 1000 elements, there are over 3 million training pairs). We then use stochastic gradient descent, with a batch size 128, to optimize the network. We use mean square error as the loss function and \textit{Adam}~\cite{kingma2014adam} as the optimizer. The implementation is done using \textit{Keras} \cite{chollet2015keras} on \textit{Tensorflow} \cite{tensorflow2015-whitepaper} and the training is conducted on a single \textit{NVIDIA Tesla} K40c graphics card.

\textbf{Guiding mesh generation via MeshingNet.}
%The training data includes a wide range of representative input parameters combined with the computed error using the error estimation approach described above. The truth solution may be computed once for multiple members of the training set, where the domain and parameters are unchanged but the given specified point varies. 
After training, the network is able to predict the target distribution $A(x)$ on a previously unseen polygonal domain. Given a problem, a LDUM is first generated; next, MeshingNet predicts the local target area upper bound, $A(x)_i$, at the centre, $x$, of the $i$th element; \textit{Triangle} then refines the LDUM to generate the non-uniform mesh based upon $KA(x)_i$ within each element of the LDUM (to be refined). Optionally, this last step may be repeated with an adjusted value of $K$ to ensure that the refined mesh has a desired total number of elements (this allows an automated approximation to ``the best possible mesh with X elements'').

\textbf{Error norms.}
Error estimation provides a means of quantifying both the local and the global error in an initial FE solution. However the precise magnitude and distribution of the error depends on the choice of the norm used to compute the difference between the LAS and the HAS. Different norms lead to different non-uniform meshes. Consequently, for any given PDE system, a single norm should be selected in order to determine $A$ from $\Gamma$, $B$ and $M$ (Equation~(\ref{mapping})). The appropriate choice of norm is a matter for the user  to decide, similar to the choice of specific \emph{a posteriori} error estimate in the conventional adaptive approach. In the following section two different norms are considered for the purposes of illustration.
%for Poisson's equation we use the L1 norm; and for the linear elasticity example we use the strain energy norm.
%This is the ``natural norm'' for the linear elasticity problem since the PDEs are the Euler-Lagrange equations for the minimization of the corresponding energy functional. Consequently, due to the linearity of the problem, the relative accuracy of two finite element solutions may be determined equivalently by which has the lower error in the energy norm or which the lower total potential energy. This property is exploited as part of our validation in subsection~\ref{Elasticity} below. 

\section{Validation results}
We now assess the performance of MeshingNet through two computational examples: a single PDE, for which we consider only the effect of the domain geometry on the optimal mesh; and a system of PDEs, for which we consider the influence of geometry, BCs and material parameters on the optimal mesh.

\subsection{2D Poisson's equation}

\begin{figure}[]
\includegraphics[width = 0.37\linewidth]{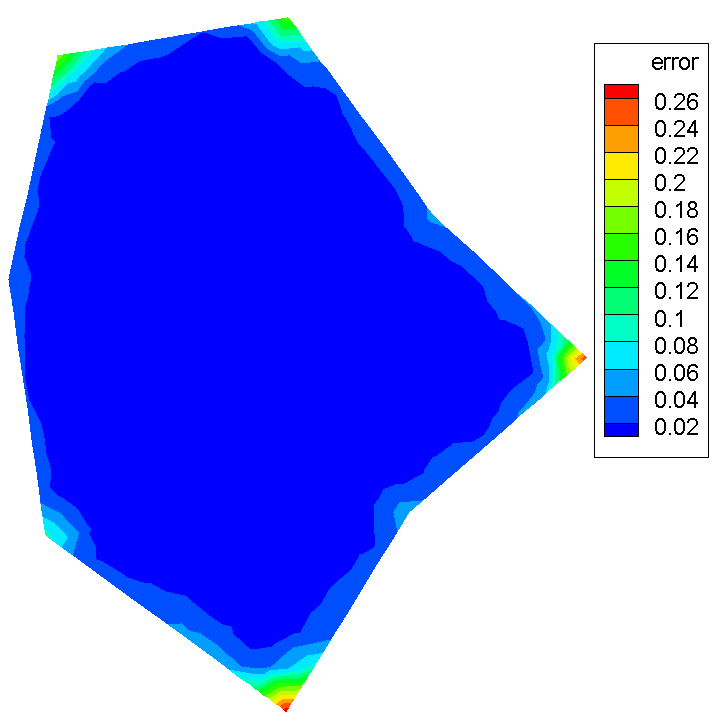}
\includegraphics[width = 0.3\linewidth]{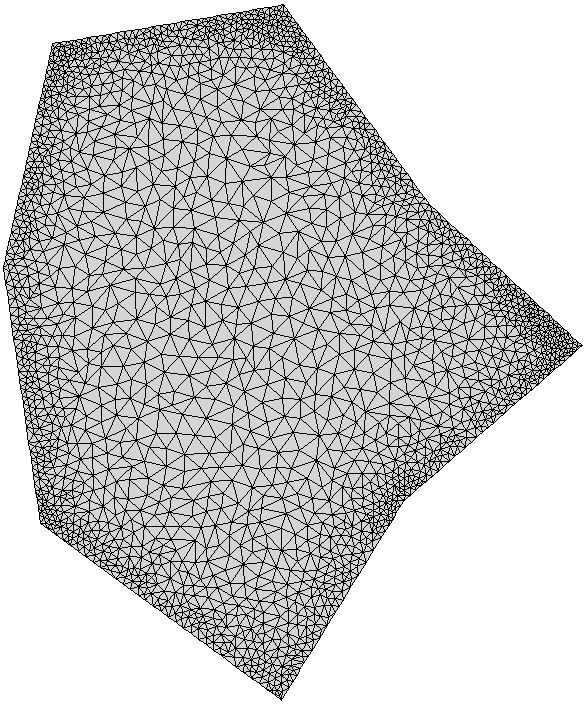}
\includegraphics[width = 0.31\linewidth]{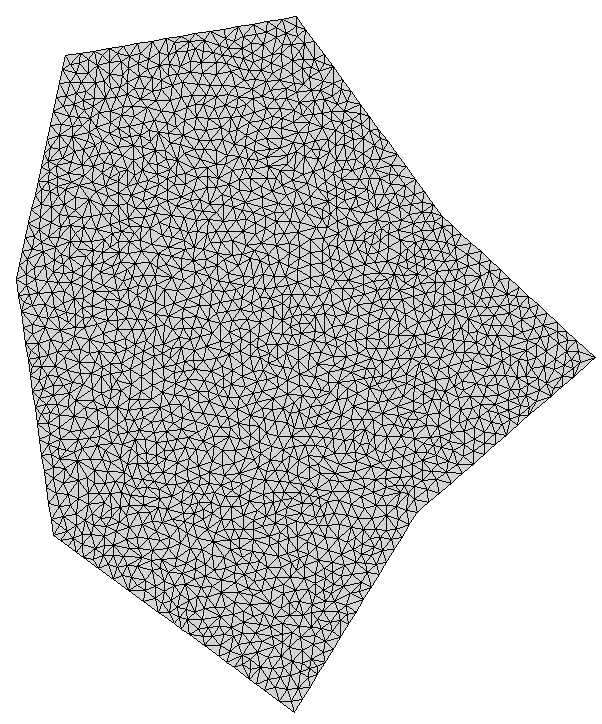}
\caption{For Poisson's equation: the L1 error distribution for the LAS (left); the 4000 elements mesh generated by \textit{Triangle} under the guidance of MeshingNet (middle); and the uniform mesh with 4000 elements generated by \textit{Triangle} (right).}
\label{fig5}
\end{figure}

\begin{figure}[]
\includegraphics[width = 0.5\linewidth]{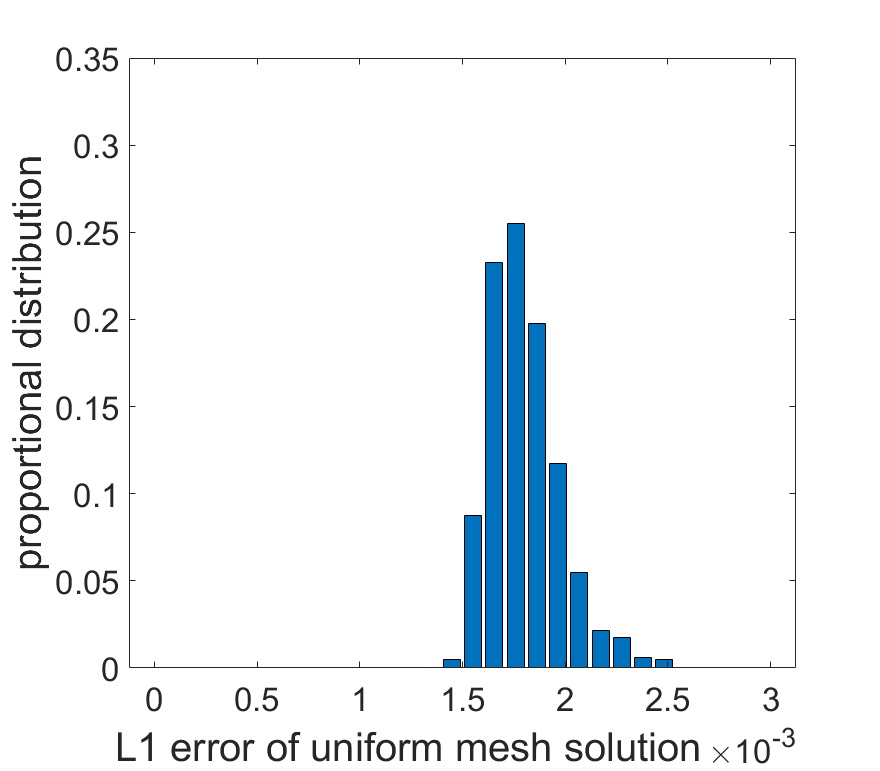}
\includegraphics[width = 0.5\linewidth]{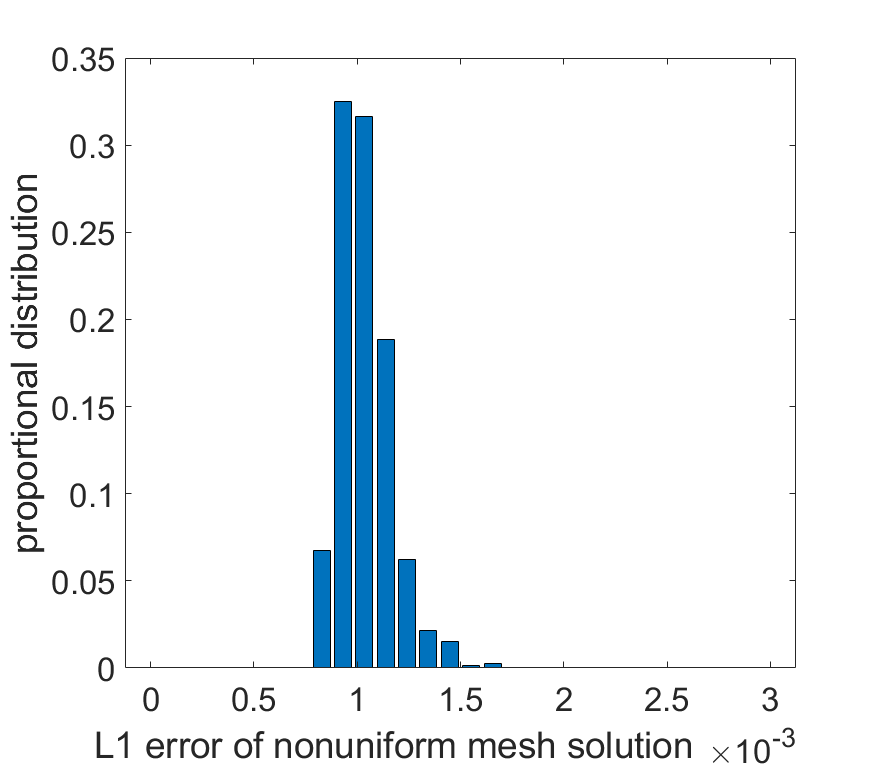}
\caption{For Poission's problem, L1 error distribution on uniform meshes (4000 elements) and MeshingNet meshes (4000 elements). Each bar shows the proportion of test data meshes whose FE solution error is in the range of the bar: the uniform meshes give FE errors between 0.0015 to 0.0025, whilst the MeshingNet meshes give FE errors between 0.0007 and 0.0015.}
\label{fig1}
\end{figure}

We solve Poisson's equation ($\nabla^2 u+1=0$) on a simply connected polygon $\Omega$ with boundary $\partial \Omega$ (on which $u=0$).
%\begin{equation}
%  \left \{ \begin{array}{lr}
%        \nabla^2 u+1=0, & x \in\Omega\\
%        u=0, & x\in\Gamma\\
%  	\end{array}\right .
%\end{equation}
The polygons in our data set are all octagons, generated randomly, subject to constraints on the polar angle between consecutive vertices and on the radius being between 100 to 200 (so the polygons are size bounded).
The L1 norm relative error estimate is
\begin{equation}
E=\left|\frac{u^{LAS}-u^{HAS}}{u^{HAS}}\right| \;.
\end{equation}
As expected, in Fig.~\ref{fig5} (which shows a typical test geometry), the mesh generated by MeshingNet is dense where the error for the LAS is high and coarse where it is low. Fig.~\ref{fig1} quantifies the improvement of MeshingNet relative to an uniform mesh (Fig.~\ref{fig5} right) by showing, for the entire test data set, the error distributions of the computed FE solutions (relative to the ground truth solutions) in each case.

\subsection{2D linear elasticity}
\label{Elasticity}
We solve 2D plane stress problems on a set of \emph{different} polygons (6-8 edges). Each polygon edge is associated with one of three possible BCs. BC1: zero displacement; BC2: uniformly distributed pressure or traction (with random amplitude up to $1000$); and BC3: unconstrained. For different geometries, we number the vertexes and edges anti-clockwise with the first vertex always on the positive x-axis, without loss of generality. To get a combinations of BCs, we always apply BC1 on the first edge, BC2 on the fourth and fifth edges and BC3 on the rest. We also allow different (homogeneous) material properties: density up to a value of $1$ and Poisson's ratio between $0$ to $0.48$. 
The error approximation uses energy norm
\begin{equation}
E=(\boldsymbol{\epsilon}^L-\boldsymbol{\epsilon}^H):(\boldsymbol{\sigma}^L-\boldsymbol{\sigma}^H)
\end{equation}
where $\boldsymbol{\epsilon}^L$ and $\boldsymbol{\epsilon}^H$ are strains of LAS* and HAS,  $\boldsymbol{\sigma}^L$ and $\boldsymbol{\sigma}^H$ are stresses of LAS* and HAS.
This is the ``natural norm'' for this problem since the PDEs are the Euler-Lagrange equations for the minimization of the following energy functional:
\begin{equation}
Ep=\int{\frac{1}{2}\boldsymbol{\epsilon}:\boldsymbol{\sigma}-\boldsymbol{F}\cdot \boldsymbol{u} d\Omega} - \int{\boldsymbol{\sigma}\cdot\boldsymbol{u}}d\Gamma
\end{equation}
where $\boldsymbol{F}$ is the body force and $\boldsymbol{u}$ is the displacement. Due to the linearity of the problem, the relative accuracy of two FE solutions may be determined equivalently by which has the lower error in the energy norm or which has the lower total potential energy (we exploit this in our validation below).

\begin{figure}[]
\includegraphics[scale =0.8]{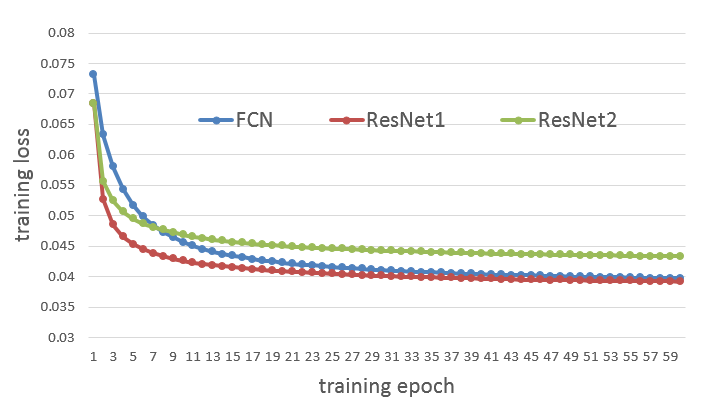}
\caption{L2 training loss on elasticity training dataset during 60 training epochs. Three curves representing FCN (blue), ResNet1 (red) and ResNet2 (green) converge individually.}
\label{fig6}
\end{figure}

%When the geometry is heptagon or hexagon we transfer them to octagon format, by locating the last (for heptagon) or last two (for hexagon) octagon vertexes on the last edge of heptagon or hexagon. There are 27 parameters as MeshingNet input.
There are 27 dimensions in MeshingNet's input: 16 for the polygon vertices, 8 for the mean value coordinates of the target point, and 1 each for the traction BC magnitude, density and Poisson's ratio. We train FCN, Resnet1 and Resnet2 for 50 epochs, each taking 142, 134 and 141 minutes respectively. Fig.~\ref{fig6} shows that the training processes all converge, with ResNet training typically converging faster than FCN.
After training, predicting the target $A$ (on all LDUM elements) for one problem takes 0.046 seconds on average, which is over 300 times faster than using the \emph{a posteriori} error method that generates the training data set.

\begin{figure}[]
\includegraphics[width = 0.336\linewidth]{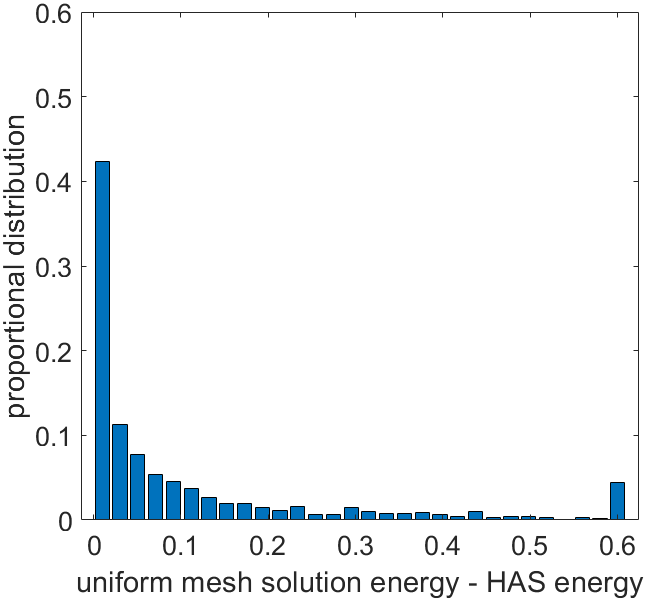}
\includegraphics[width = 0.333\linewidth]{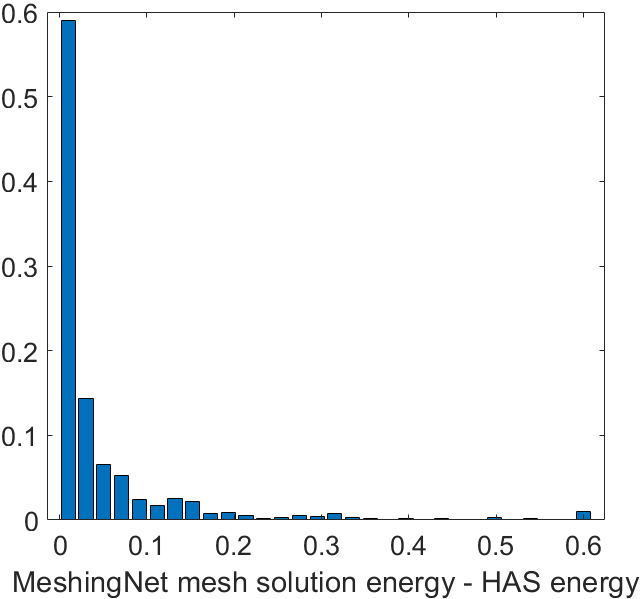}
\includegraphics[width = 0.31818\linewidth]{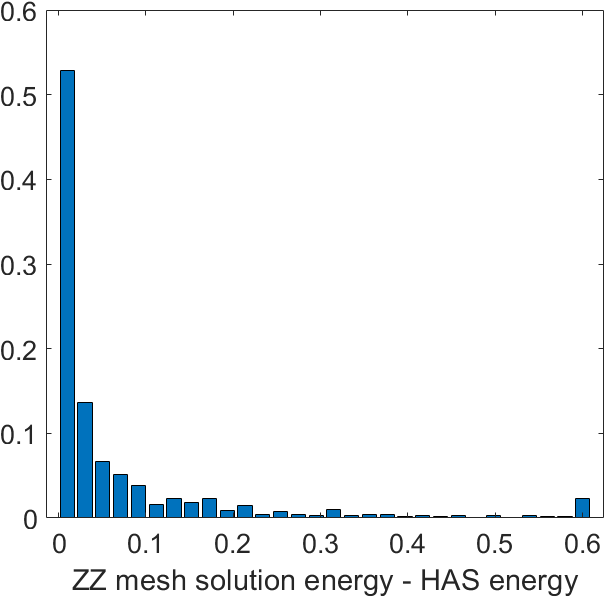}
\caption{Potential energy comparison of solving 2D elasticity test problems on 4000 element meshes: uniform element size (left); MeshingNet (using FCN) meshes (centre); and ZZ refined meshes (right). We use the HAS energy as our baseline: the MeshingNet mesh solutions have energies that are significantly closer to the HAS energies than the ZZ mesh solutions and uniform mesh solutions (since a greater proportion of results are distributed near zero). The rightmost bar represents the proportion of all tests where the energy difference is no smaller than 0.6.}
\label{fig2}
\end{figure}

\begin{figure}[!]
%Heptagon Traction: -186			Poission's ratio: 0.0717			density: 0.2536
%\includegraphics[width = 0.5\linewidth]{3001e.png}
%\includegraphics[width = 0.5\linewidth]{3001m.png}
 
\includegraphics[width = 0.32\linewidth]{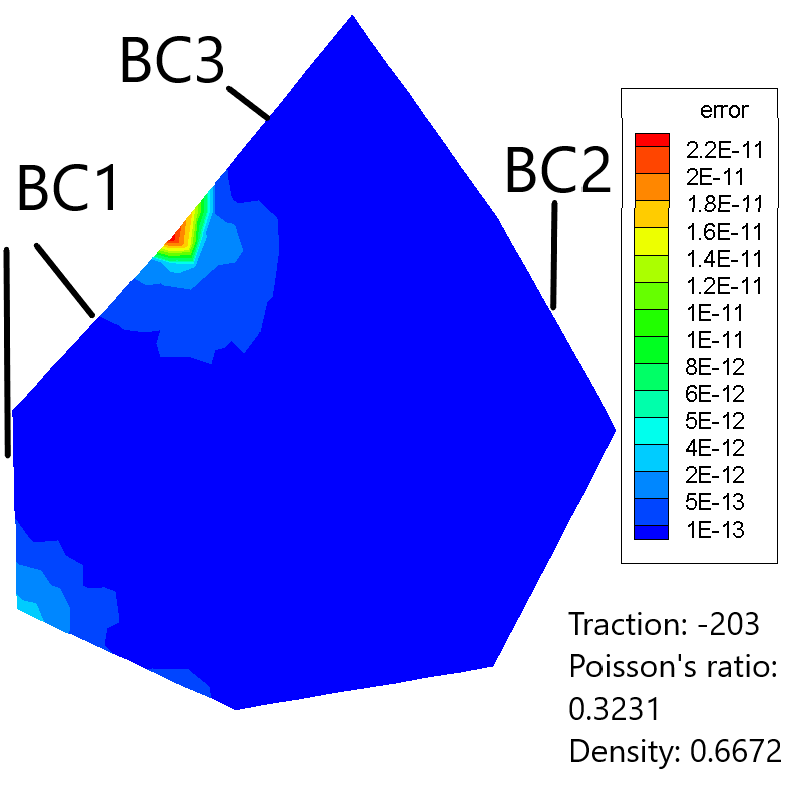}
\includegraphics[width = 0.32\linewidth]{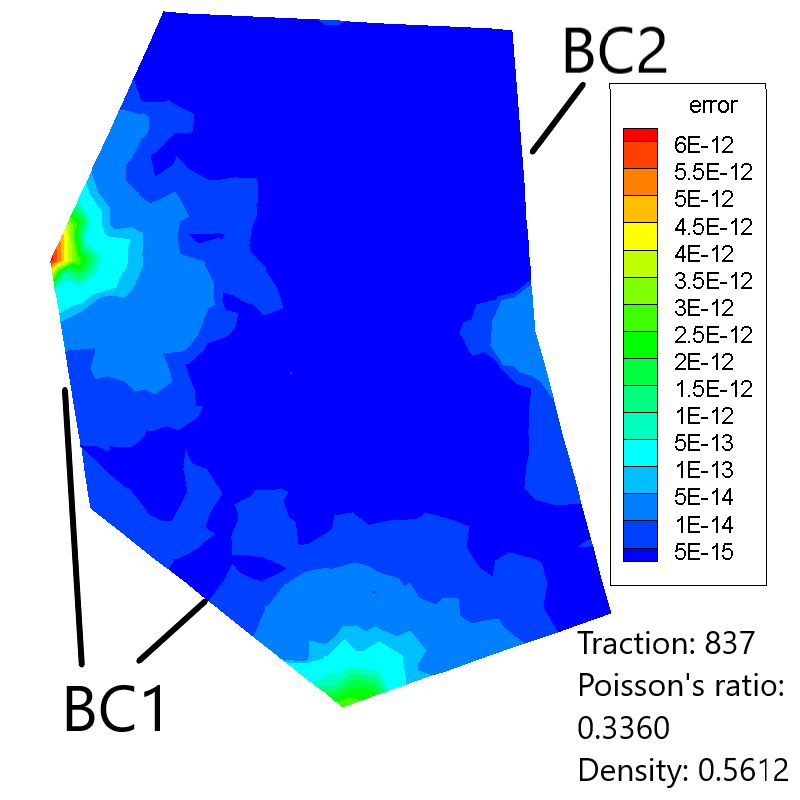}
\includegraphics[width = 0.33\linewidth]{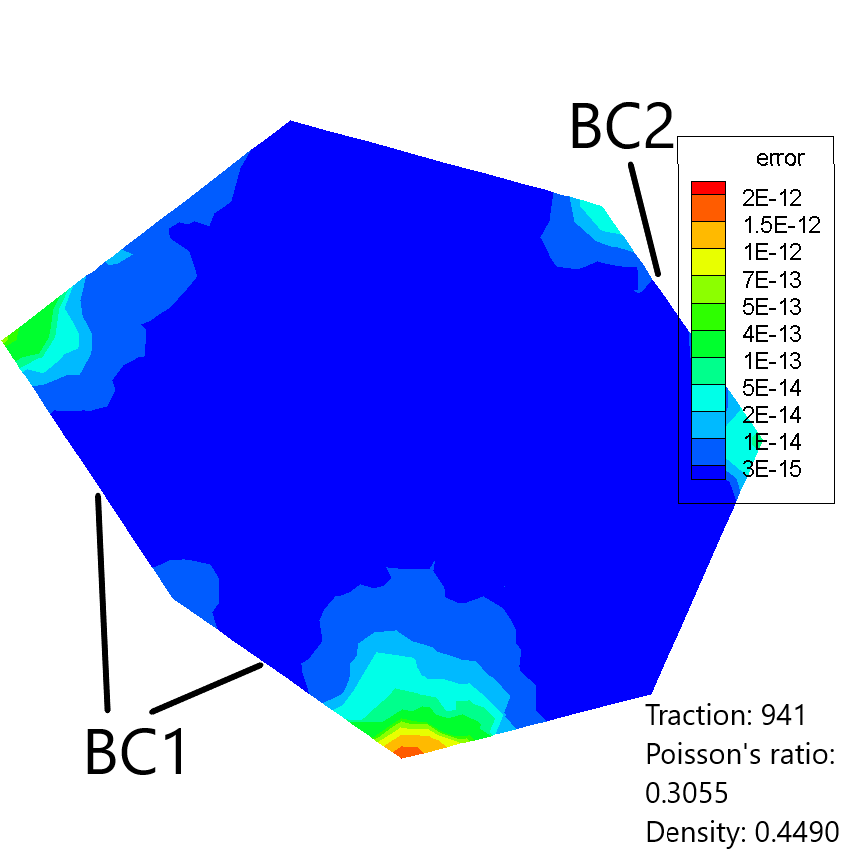}

\includegraphics[width = 0.32\linewidth]{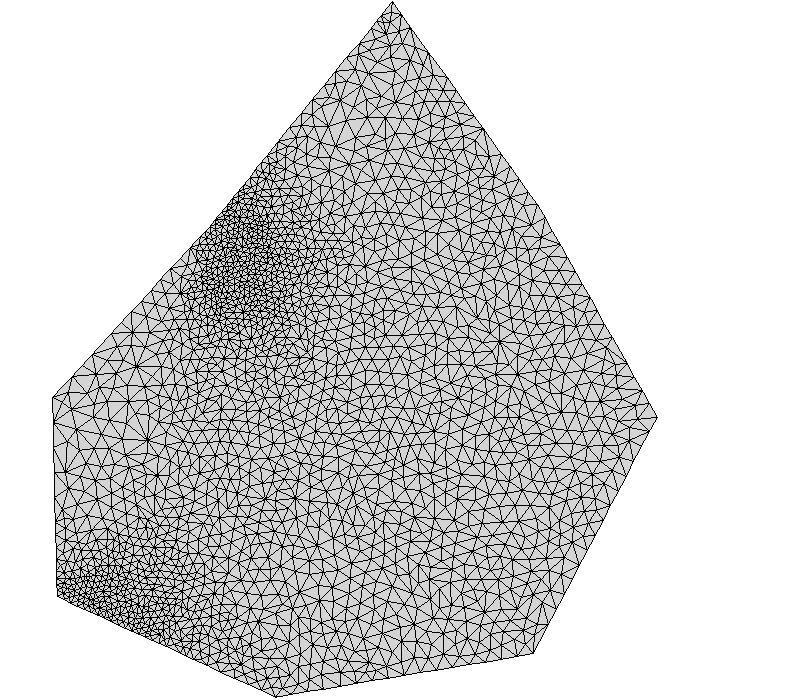}
\includegraphics[width = 0.32\linewidth]{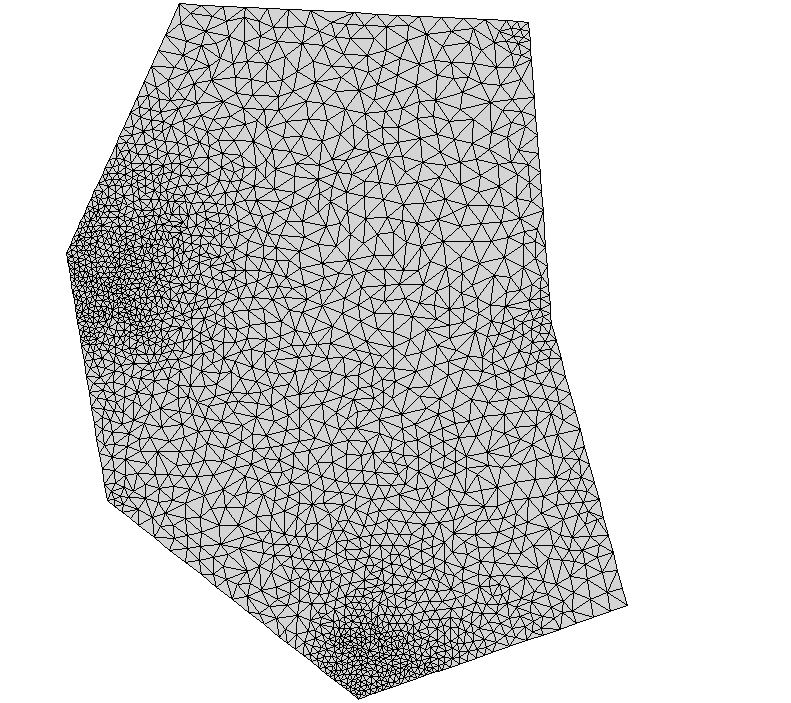}
\includegraphics[width = 0.33\linewidth]{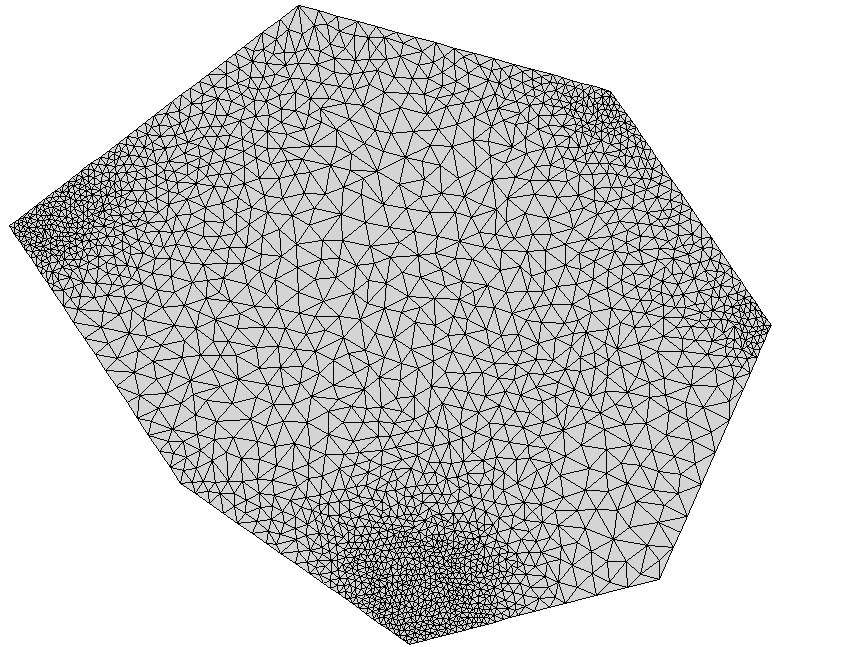}

\includegraphics[width = 0.32\linewidth]{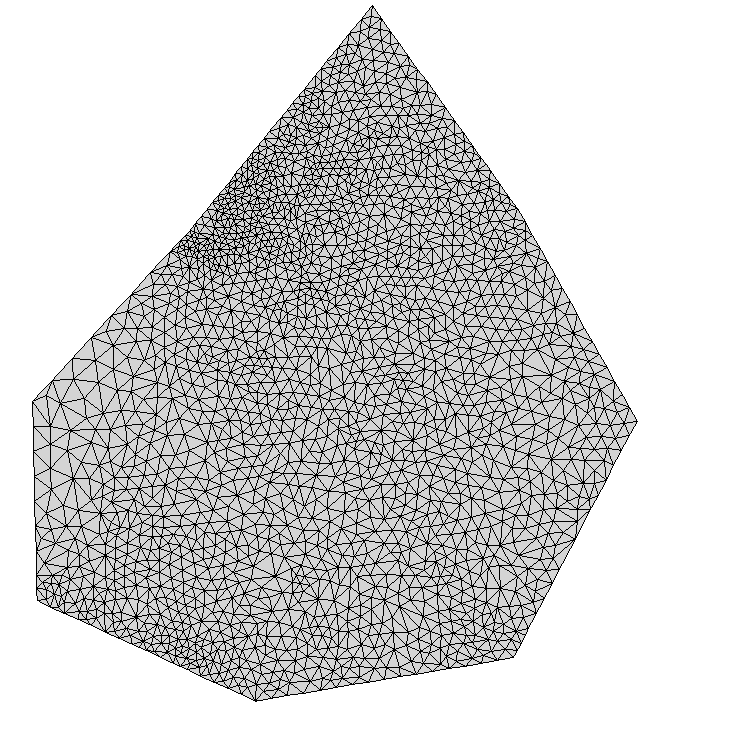}
\includegraphics[width = 0.32\linewidth]{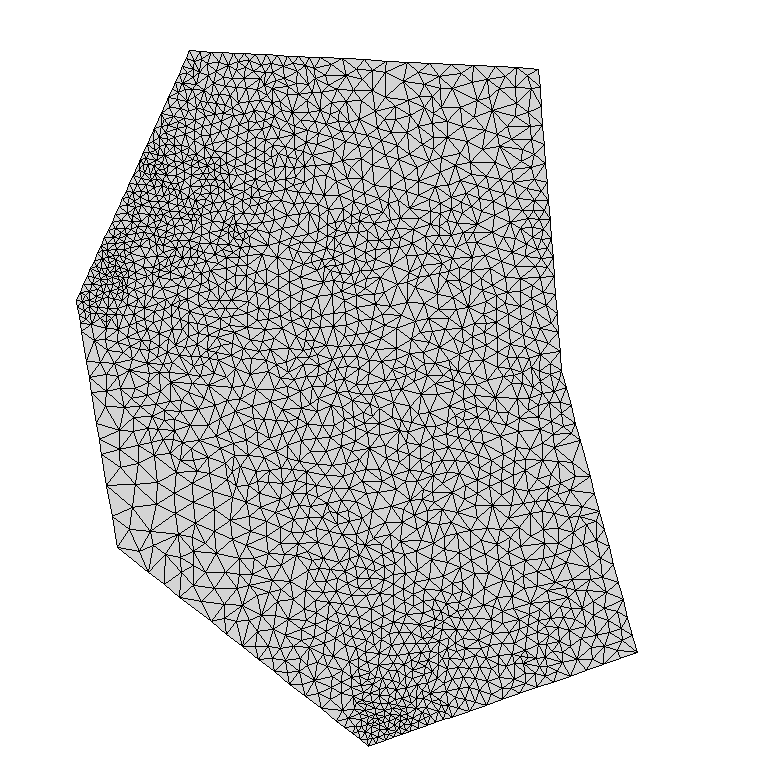}
\includegraphics[width = 0.33\linewidth]{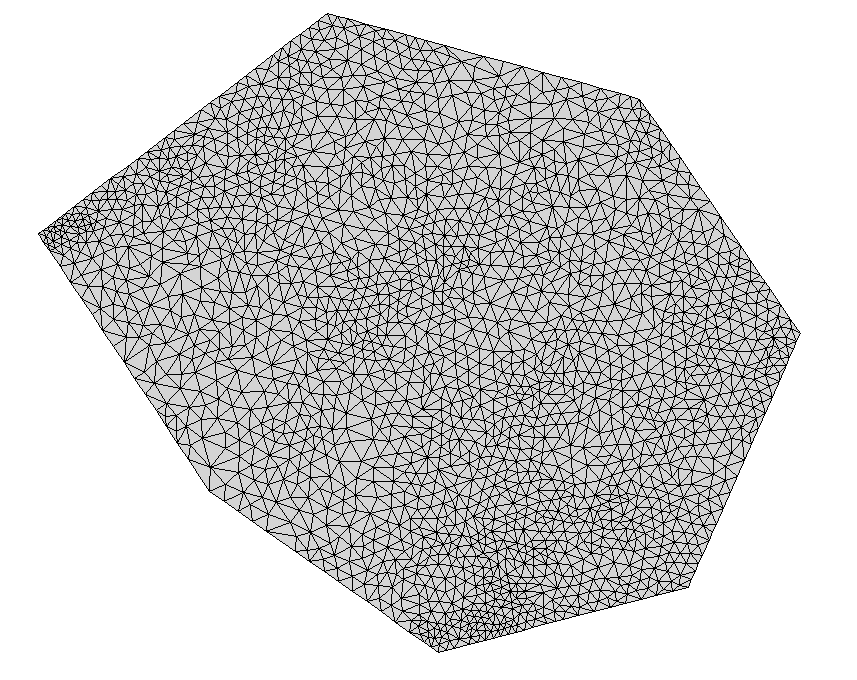}

\caption{FE error (top) relative to HAS on coarse uniform mesh, non-uniform mesh guided by MeshingNet (middle) and non-uniform mesh refined by ZZ (bottom). The left geometry is an octagon and other two are heptagons (defined by placing their final vertex at the centre of edge 7).}
\label{fig3}
\end{figure}

Fig.~\ref{fig2} shows a comparison of FE results computed on MeshingNet meshes, uniform meshes of the same number of elements (4000 elements) and non-uniform meshes (also of the same number of elements) computed based upon local refinement following ZZ error estimation. The former meshes have FE solutions with potential energy significantly lower than the uniform mesh and ZZ refined mesh (and much closer to the HAS potential energy).
Fig.~\ref{fig3} illustrates some typical meshes obtained using MeshingNet: the non-uniform meshes correspond to the error distributions in the LAS. Though not shown here due to space constraints, we also find that the traction-to-density ratio impacts the non-uniform mesh most significantly. Overall, this example shows that MeshingNet can generate high quality meshes that not only account for geometry but also the given material properties and BCs.
%This experiment demonstrates that MeshingNet aids generating efficacious non-uniform mesh on unseen geometries.


Finally we compare different DNN models, using the average potential energy on the testing data set. The baseline is from the HDUMs, whose FE solutions have a mean energy of $-7.7293$, followed by the meshes from ResNet2 (mean energy $-7.6816$), ResNet1 ($-7.6815$), and FCN ($-7.6813$). The lower the better. These are all far superior to the uniform meshes of the same size (4000 elements), which yield FE solutions with a mean energy of $-7.6030$. Note that ResNet not only shortens the training time over FCN but, on average, produces better solutions.
%Table~\ref{tab1} gives the average potential energy, indicates the neural network structure of MeshingNet has effect to the quality of the generated mesh. We find ResNet performs lower potential energy than FCN even when FCN has more parameters to train. With less parameters ResNet can shorten the training time while giving higher quality mesh than FCN.

%\begin{table}
%\caption{Average potential energy on 800 testing data}\label{tab1}
%\begin{tabular}{ |p{6cm}|p{6cm}|  }%{|l|l|}
%\hline
%Meshes & average potential energy\\
%\hline
%FCN nonuniform mesh &  -7.6813\\
%ResNet1 nonuniform mesh &  -7.6815\\
%ResNet2 nonuniform mesh & -7.6816 \\
%uniform mesh & -7.6030\\
%HDUM & -7.7293\\
%\hline
%\end{tabular}
%\end{table}

\section{Discussion}
In this paper, we have proposed a new non-uniform mesh generation method based on DNN. The approach is designed for general PDEs with a range of geometries, BCs and problem parameters. We have implemented a two-dimensional prototype and validated it on two test problems: Poisson's equation and linear elasticity. These tests have shown the potential of the technique to successfully learn the impact of domain geometries, BCs and material properties on the optimal finite element mesh. Quantitatively, meshes generated by MeshingNet are shown to be more accurate than uniform meshes and non-uniform ZZ meshes of the number of elements. Most significantly, MeshingNet avoids the expense of \emph{a posteriori} error estimation whilst still predicting these errors efficiently. Even though generating the training data set is expensive, it is offline and is thus acceptable in practice.

The meshes generated via MeshingNet may be used in a variety of ways. If our goal is to obtain a high quality mesh with a desired number of elements then the approach described in this paper provides a cost-effective means of achieving this. If however the goal is to produce a solution with an estimated {\em a posteriori} error that is smaller than a desired tolerance, then the generated mesh may not meet this criterion. This may be addressed either by regenerating the mesh based upon a higher target number of elements in the final mesh, or through the use of a traditional {\em a posteriori} estimate on the computed solution in order to guide further mesh refinement. In the latter case we can view MeshingNet as a means to obtaining an improved initial mesh within a traditional mesh adaptivity workflow.

%We are expecting to both extend the complexity of PDEs and increase utility of MeshingNet. We attempt to generate nonuniform mesh for 3D elasticity problem by taking after methodology in 2D elasticity. Besides, we have interest in 2D and 3D, steady and transient fluid flow problems. For utility, we consider an automatic program that transfer computer aided design format geometry into “feed" of neural networks. Prospectively, mesh can be generated for interactions of elastic assemblies or multi phase objects. The greatest challenge for these problems is how to enhance the approximation capability and accuracy of neural networks.  Generally, people use convolutional neural networks for high dimensional input and recurrent network for time dependant problems. However, cautions should be taken that huge network structures usually increase both training and running time.

In future, we plan to generalize MeshingNet onto more general problems: 3D geometries, more complex PDE systems and BCs (e.g. Navier-Stokes and fluid-structure interactions), and time-dependent cases.
For 3D problems, \textit{Tetgen}~\cite{si2015tetgen} is able to do a similar refinement process to what \textit{Triangle} does in 2D in this paper, which will enable us to directly apply MeshingNet to 3D problems. Furthermore, it would be desirable to develop an interface to enable the mesh generator to read geometries from standard computer aided design software. For complex three-dimensional problems it also seems unlikely that the accurate, but very expensive, approach to training the error estimator that is used in this paper will always be computationally viable (despite its excellent performance). In such cases we may replace this estimator with a more traditional method such as~\cite{ainsworth1997posteriori,bank1985some,zienkiewicz1991adaptivity}, on relatively fine uniform grids for training purposes.

Finally, we note that the ANNs used in this initial investigation are relatively simple in their structure. In the future the use of new DNN models, such as Convolutional/Graph Neural Networks, should be considered. These may be appropriate for problems in three dimensions or with larger data sets, such as arising from more general geometries and boundary conditions.


%\clearpage
\bibliographystyle{splncs04}
\bibliography{reference}

%\end{thebibliography}
\end{document}